\documentclass {amsart}

\usepackage {amssymb, amsmath, amsthm}

\newtheorem {Te}[equation]{Theorem}
\newtheorem {Le}[equation]{Lemma}
\newtheorem {Ps}[equation]{Proposition}

\theoremstyle {definition}
\newtheorem {pag}[equation]{}

\makeatletter \@addtoreset{equation}{section} \makeatother

\title {Bound for the order for $p$-elementary subgroups in the plane Cremona group over a perfect field}
\author {A.~L.~Fomin}
\thanks {The author was partially supported by Leading Scientific Schools, grant No.~4713.2010.1}
\address {Andrei Fomin: Department of Higher Algebra, Faculty of Mechanics and Mathematics, Moscow State Lomonosov University, Vorobievy Gory, Moscow, 119~899, Russia}
\email {fin-al@yandex.ru}

\begin {document}

\maketitle

\begin {abstract}
We obtain a sharp bound for $p$-elementary subgroups in the Cremona group $\mathrm {Cr}_2 (k)$ over an arbitrary perfect field $k$.
\end {abstract}

\section {Introduction}

Let $k$ be a field. The Cremona group $\mathrm {Cr}_n (k)$ over $k$ is the group of birational transformations of $\mathbb P^n$ that are defined over $k$. Though studying finite subgroups of $\mathrm {Cr}_2 (\mathbb C)$ has a long history that has lasted for nearly a century and a half dealing with fields $k$, which are not algebraically closed, started only a few years ago in~\cite {DI}.  

A finite abelian group $A$ is called a \emph {$p$-elementary group}, where $p$ is a prime number, if $A \cong (\mathbb Z/p)^r$, $r$ is called the rank of $A$ and denoted by $\mathrm {rk}\,A$. In~\cite {Be} A.~Beauville classified maximal $p$-elementary subgroups in $\mathrm {Cr}_2 (k)$ over an algebraically closed field $k$ of arbitrary characteristic up to conjugacy. The purpose of the present paper is to find a sharp bound for $p$-elementary subgroups in the plane Cremona group $\mathrm {Cr}_2 (k)$ over an arbitrary perfect field $k$.

For a perfect field $k$, denote by $\overline k$ its algebraic closure and put $\Gamma_k = \mathrm {Gal} (\overline k/k)$. For a prime number, $p$ it is always assumed that $p \ne \mathrm {char}(k)$. Define $t = [k(\zeta_p):k]$, where $\zeta_p \in \overline k$ is any primitive root of unity of degree $p$. It is clear that $t$ divides $p-1$.

Our main result is the following

\begin {Te} \label{main-theorem}
Let $G \subset \mathrm {Cr}_2 (k)$ be a $p$-elementary subgroup, where $k$ is a perfect field. Then

\begin{equation} \label{equation}
\mathrm {rk}\,A \le
\begin {cases}
4 & \text {if $p=2$;}\\
3 & \text {if $p=3,\, t=1$;}\\
2 & \text {if $p=3,\, t=2$ and $p>3,\, t=1,2$;}\\
1 & \text {if $t=3,4,6$;}\\
0 & \text {otherwise.}\\ 
\end {cases}
\end{equation}
Moreover, this bound is attained for any $p \ne \mathrm {char}(k)$.
\end {Te}

\section {Bounds for a $p$-torsion subgroup of a torus}
\begin{pag}
Let $T$ be an algebraic torus of dimension $d$ defined over $k$. In~\cite {Se} J.-P.~Serre obtained a sharp bound for the order of finite $p$-subgroups in $T(k)$. Below we give a similar bound for $p$-elementary subgroups.
\end{pag}

\begin {Te}
\label{theorem2}
In the notation above, $\mathrm {rk}\,T(k)[p] \le \frac d {\varphi (t)}$, where $\varphi$ is Euler's function. Moreover, this bound is attained for a suitable torus defined over $k$.
\end {Te}

\begin {proof}
Let $\mathrm X (T)$ and $\Upsilon (T)$ be the groups of characters and cocharacters of $T$ over $\overline k$, $\rho \colon \Gamma_k \to \mathrm {Aut} (\Upsilon (T))$ the action of the Galois group, $\rho_p \colon \Gamma_k \to \mathrm {Aut} (\Upsilon (T)/p)$ its reduction modulo $p$. 
Let also $\boldsymbol \mu_p \subset \overline k^*$ be the group of the roots of unity of degree $p$, and $\chi \colon \Gamma_k \to \mathrm {Aut} (\boldsymbol \mu_p) \cong (\mathbb Z/p)^*$ the action of the Galois group.

It is clear that $T(k)[p] = T(\overline k)[p]^{\Gamma_k}$ and $T(\overline k)[p] \cong \mathrm {Hom}\,(\mathrm X (T)/p,\; \boldsymbol \mu_p) \cong \Upsilon (T)/p \otimes \boldsymbol \mu_p$,
all isomorphisms being compatible with the actions of the Galois group. Obviously, $\mathrm {rk}\,(\Upsilon (T)/p \otimes \boldsymbol \mu_p)^{\Gamma_k} \le \mathrm {rk}\,(\Upsilon (T)/p \otimes \boldsymbol \mu_p)^g$ for any $g \in \Gamma_k$ and $g$ acts on $\Upsilon (T)/p \otimes \boldsymbol \mu_p$ as $\rho_p (g) \otimes \chi (g) = \chi (g) \rho_p (g) \otimes 1$. Using any isomorphism $\boldsymbol \mu_p \cong \mathbb Z/p$ and $\Upsilon (T)/p \otimes \boldsymbol \mu_p \cong \Upsilon (T)/p$ it is possible to identify the set of fixed points of $g$ in $\Upsilon (T)/p \otimes \boldsymbol \mu_p$ with the set of fixed points of $\chi (g) \rho_p (g)$ in $\Upsilon (T)/p$, which is nothing but the eigenspace of $\rho_p (g)$ with eigenvalue $\chi (g)^{-1}$.

We fix $g \in \Gamma_k$ such that $\chi (g)$ is of order $t$ and denote $\chi (g)^{-1} = \varepsilon$. Since $\rho (g)$ has finite order, its characteristic polynomial $F$ is the product of cyclotomic polynomials: $F = \prod\limits_i \Phi_{d_i}$, and the characteristic polynomial of $\rho_p (g)$ is $\overline F = \prod\limits_i \overline \Phi_{d_i}$, where $\overline \Phi$ denotes the reduction a polynomial $\Phi$ modulo $p$. To prove the theorem, we need to find an upper bound for the multiplicity of $\varepsilon$ as the root of $\overline \Phi_{d_i}$.

\begin {Le}
\label{lemma}
In the above notation, the multiplicity of $\varepsilon \in (\mathbb Z/p)^*$ as the root of $\overline \Phi_n$ is the same for all $\varepsilon$ of the fixed order $t$, and it is positive iff $n = tp^f$. 
\end {Le}

\begin {proof}[Proof of the Lemma] 
First, if $p \nmid n$ and $q = p^f$, then $\overline \Phi_{nq} \equiv \overline \Phi_n^{\varphi (q)} \pmod p$, so we can assume that $p \nmid n$.

Let $\mathcal O$ be the integral closure of $\mathbb Z$ in the field $\mathbb Q (\zeta_n)$, where $\zeta_n \in \mathbb C$ is any primitive root of unity of degree $n$, $\boldsymbol \mu_n \subset \mathcal O^*$ the group of the roots of unity of degree $n$, and $\mathfrak p \subset \mathcal O$ any prime ideal such that $\mathfrak p \cap \mathbb Z = p \mathbb Z$. Then $\Phi_n (X) = \prod\limits_{\zeta} (X - \zeta)$ and $\overline \Phi_n (X) = \prod\limits_{\zeta} (X - \overline \zeta)$ in $\mathcal O/\mathfrak p$, where $\zeta$ runs through all primitive roots of unity of degree $n$. 
It is well known that the natural map $\boldsymbol \mu_n \to (\mathcal O/\mathfrak p)^*$ is injective, so $\overline \zeta$ is of order $n$ in $(\mathcal O/\mathfrak p)^*$ for any $\zeta$. This implies that the set of roots of $\overline \Phi_n$ in $\mathcal O/\mathfrak p$ coincides with the set of all elements of order $n$ in $(\mathcal O/\mathfrak p)^*$.

Suppose that $\overline \Phi_n$ has a root $\varepsilon \in (\mathbb Z/p)^*$ of order $t$; then $t=n$ and any element of order $t$ in $(\mathbb Z/p)^*$ is a simple root of $\overline \Phi_n$. 
This proves all statements of the lemma.    
\end {proof}

Going back to the proof of Theorem~\ref{theorem2} we see that it follows from the above lemma 
that the multiplicity of $\varepsilon$ as the root of $\overline \Phi_{d_i}$ is bounded from above by $\frac {\varphi (d_i)}{\varphi (t)}$, and its multiplicity as the root of $\overline F$ is bounded from above by $\frac d{\varphi (t)}$, since $\sum\limits_i \varphi (d_i) = d$.

To prove the second statement of Theorem~\ref{theorem2}, it is enough to construct a torus of dimention $d = \varphi (t)$ defined over $k$ such that $\mathrm {rk}\,T(k)[p] > 0$. This is done in~\cite {Se} (see the proof of Theorem 4$'$).                
\end {proof}
 
\section {Proof of the main theorem}

In this section we prove Theorem~\ref{main-theorem}.
\begin{pag} \label{division-into-cases}
Let $A \subset \mathrm {Cr}_2 (k)$ be a $p$-elementary subgroup. It is known (see~\cite [Theorem~5]{DI}) that $A$ can be represented as a subgroup of $\mathrm {Aut}_k (S)$, where $S$ is a smooth projective surface defined and rational over $k$, which belongs to one of the following two types:

\begin {enumerate}
\item There exists an $A$-equivariant conic bundle structure $f \colon S \to \mathbb P^1$ such that $\mathrm {rk\,Pic} (S/\mathbb P^1)^A = 1$; 

\item $S$ is a Del Pezzo surface such that $\mathrm {rk\,Pic} (S)^A = 1$.
\end {enumerate}    
\end{pag}

\begin {Ps} \label{proposition}
If $p \nmid n$, any $p$-elementary subgroup $A \subset G(k)$, where $G$ is a $k$-form of $\mathbf {PGL}_n$, is contained in a maximal torus defined over $k$. 
\end {Ps}

\begin {proof}
This statement was proved in~\cite [Lemma~3.1]{Be} for $k = \overline k$. The centralizer of $A$ in $G$, which is defined over $k$ as $A$ itself is, contains a maximal torus defined over $k$, which is the maximal torus in $G$. Since $A$ consists of semisimple elements, any maximal torus that centralizes $A$ must contain it. 
\end {proof}

\begin{pag}
In the case $\mathrm {rk}\,A \ge 1$, it was proved in~\cite [Theorem~2]{DI} that $t \in \{ 1,\,2,\,3,\,4,\,6 \}$ and moreover for these values of $t$ the desired $A$ does exist.
\end{pag}

\begin{pag}
Suppose that $\mathrm {rk}\,A \ge 2$. We shall prove that $t \le 2$. We can assume that $p>3$ as otherwise there is nothing to prove and that $A$ is a subgroup of $\mathrm {Aut}_k (S)$ as it is described above. Denote $\overline S = S \otimes \overline k$. We have two cases of~\ref{division-into-cases}.

Let $f \colon S \to \mathbb P^1$ be an $A$-equivariant conic bundle. The action of $A$ on the base defines the homomorphism $A \to \mathrm {Aut}_k (\mathbb P^1)$. Denote by $\overline A$ its image and by $A_0$ its kernel. Obviously, $A_0$ is an automorphism group of the generic fiber of $f$, which is a smooth curve of genus $0$ over the field $K$ of rational functions on $\mathbb P^1$. The automorphism group of the base is isomorphic to $\mathbf {PGL}_2$ over $k$, and the automorphism group of the generic fiber is a $K$-form of $\mathbf {PGL}_2$. It is readily seen that $t$ has the same value for $k$ and $K$. Since $p$ is odd, it follows from Proposition~\ref{proposition} that $\overline A$ and $A_0$ are contained in tori of dimension $1$ defined over $k$ and $K$ respectively. Theorem~\ref{theorem2} yields that $\mathrm {rk}\,A_0 \le 1$ and $\mathrm {rk}\,\overline A \le 1$, with the equality being possible only if $t \le 2$. Finally we obtain that $\mathrm {rk}\,A \le 2$, where the equality holds only if $t \le 2$. 

Let $S$ be a Del Pezzo surface. It follows from~\cite [Proposition~3.9]{Be} and \cite [Theorem~5]{DI} that $9\ge K_S^2 \ge 6$ and $K_S^2 \neq 7$. Consider possibilities for $K_S^2$ case by case.

\begin{itemize}
 \item 
If $K_S^2 = 9$, then $\overline S \cong \mathbb P^2$ and thus $S \cong \mathbb P^2$ over $k$, 
since $S$ is rational over $k$. Therefore $\mathrm {Aut}(S) \cong \mathbf {PGL}_3$ 
over $k$ and Proposition~\ref{proposition} gives that $A$ is contained in a torus of dimension $2$. According to Theorem~\ref{theorem2} this is possible only if $t \le 2$.

 \item
If $K_S^2 = 8$, then $\overline S \cong \mathbb P^1 \times \mathbb P^1$. (Otherwise $\overline S$ contains a unique $(-1)$-curve which must be defined over $k$. This contradicts $\mathrm {rk\,Pic} (S)^A = 1$). Then the connected component $\mathrm {Aut}(S)^{\circ}$ is a $k$-form of $\mathbf {PGL}_2 \times \mathbf {PGL}_2$ of index $2$ in $\mathrm {Aut}(S)^{\circ}$. It is clear that $A \subset \mathrm {Aut}(S)^{\circ}$, and by Proposition~\ref{proposition} $A$ is contained in a torus of dimension $2$, and thus $t \le 2$.     

 \item
If $K_S^2 = 6$, then the connected component $\mathrm {Aut}(S)^{\circ}$ is a $2$-dimensional torus of index $12$ in $\mathrm {Aut}(S)$. Since $A \subset \mathrm {Aut}(S)^{\circ}$, we obtain that $t \le 2$. 
\end{itemize}

To prove the second statement of Theorem~\ref{main-theorem} for the case $t \le 2$, applying Theorem~\ref{theorem2} we obtain that for such $t$ there exists a $2$-dimensional torus $T$ defined over $k$ such that $T(k)$ contains a $p$-elementary subgroup $A$ of rank $2$. Thus the well known fact that $T$ is rational over $k$ (see~\cite [\S 4.9]{V}) yields that $A \subset \mathrm {Cr}_2 (k)$.
\end{pag}

\begin{pag}
Suppose now that $p$ is odd and $\mathrm {rk}\,A \ge 3$. It is shown in~\cite [Proposition~2.6 and Proposition~3.10]{Be} that $p=3$, $\mathrm {rk}\,A = 3$ and $S$ must be a cubic surface in $\mathbb P^3$. We claim that $t=1$.

It follows from Proposition~\ref{proposition} that $A \subset T(k)$, where $T \subset \mathbf {PGL}_4$ is a maximal torus defined over $k$. We use notation from the proof of Theorem~\ref{theorem2}. Since $\mathbf {PGL}_4$ is a group of inner type, for any $g \in \Gamma_k$ $\rho (g)$ acts on $\Upsilon (T)$ as an element of the Weyl group. Let $F = \prod\limits_i \Phi_{d_i}$ be the characteristic polynomial of $\rho (g)$ and $\overline F = \prod\limits_i \overline \Phi_{d_i}$ its reduction modulo $3$. Note that each $d_i$ divides one of the invariant degrees of the Weyl group, therefore each $d_i \in \{1,\,2,\,3,\,4\}$. Suppose that $t=2$; then the multiplicity of $-1 \in (\mathbb Z/3)^*$ as the root of $\overline F$ is equal to $3$. It follows easily from Lemma~\ref{lemma} that each $d_i = 2$ and $F(X) = (X+1)^3$. Since $\rho (g)$ has finite order, $\rho (g) = -1$, but it is well known that $-1$ does not belong to the Weyl group of $\mathbf {PGL}_4$. So we conclude that the case $t=2$ is impossible. This completes the proof of~\eqref {equation} for $p>2$.

As for the second statement of Theorem~\ref{main-theorem}, if $p=3$ and $t=1$, i.e. $k$ contains the primitive cubic root of unity, the Fermat cubic given by equation $X_0^3+X_1^3+X_2^3+X_3^3=0$ in $\mathbb P^3$ is rational over $k$ and evidently admits the action of $3$-elementary group $A$ with $\mathrm {rk}\,A = 3$, so $A \subset \mathrm {Cr}_2 (k)$.
\end{pag}

\begin{pag}
Finally, suppose that $p=2$. It was proved in~\cite [Proposition~2.6 and Proposition~3.11]{Be} that $\mathrm {rk}\,A \le 4$. On the other hand, $\mathbb P^1$ admits $(\mathbb Z/2)^2$ as the automorphism group for every field $k$, hence there exists an action of the group $A \cong (\mathbb Z/2)^4$ on $\mathbb P^1 \times \mathbb P^1$ and $A \subset \mathrm {Cr}_2 (k)$. This completes the proof of the main theorem.           
\end{pag}

\end {document}